\documentclass{siamart190516}

\usepackage{low_rank_markov_chains_style}

\headers{Low-rank tensor methods for Markov chains}{Georg, Grasedyck, Klever, Schill, Spang, Wettig}

\title{Low-rank tensor methods for Markov chains with applications to tumor progression models%
			\thanks{
				Submitted to the editors.
				\funding{This work was partially supported by the German Research Foundation (DFG)  grant~SFB/TRR-55 ``Hadron Physics from Lattice QCD'',~SPP-1886 ``Polymorphic uncertainty modelling for the numerical design of structures'', and~FOR-2127 ``Selection and Adaptation during Metastatic Cancer Progression''.
				}
				}
			}

\author{
Peter Georg\thanks{Department of Physics, University of Regensburg, 93040 Regensburg, Germany (\email{peter.georg@ur.de}, \email{tilo.wettig@ur.de}).}
\and Lars Grasedyck\thanks{Institut für Geometrie und Praktische Mathematik, RWTH Aachen University, 52062 Aachen, Germany (\email{lgr@igpm.rwth-aachen.de}, \email{klever@igpm.rwth-aachen.de}).}
\and Maren Klever\footnotemark[3] \footnotemark[2]
\and Rudolf Schill\thanks{Department of Statistical Bioinformatics, Institute of Functional Genomics, University of Regensburg,  93040 Regensburg, Germany (\email{rudolf.schill@ur.de}, \email{rainer.spang@ur.de}).}
\and Rainer Spang\footnotemark[4]
\and Tilo Wettig\footnotemark[2]
}

\begin{document}

\maketitle 

\begin{abstract}
Continuous-time Markov chains describing interacting processes exhibit a state space that grows exponentially in the number of processes.
This state-space explosion renders the computation or storage of the time-marginal distribution, which is defined as the solution of a certain linear system, infeasible using classical methods.
We consider Markov chains whose transition rates are separable functions, which allows for an efficient low-rank tensor representation of the operator of this linear system.
Typically, the right-hand side also has low-rank structure, and thus we can reduce the cost for computation and storage from exponential to linear.
Previously known iterative methods also allow for low-rank approximations of the solution but are unable to guarantee that its entries sum up to one as required for a probability distribution.
We derive a convergent iterative method using low-rank formats satisfying this condition.
We also perform numerical experiments illustrating that the marginal distribution is well approximated with low rank.
\end{abstract}



\section{Introduction}
\label{section:introduction}

The dynamics of cancer can be studied using tumor progression models, cf.~\cite{Beerenwinkel2014}.
These models describe the evolving genotype of a tumor as a continuous-time Markov chain.
A model includes $d$ genomic loci that may or may not be mutated.
It starts out with all mutations absent and then progressively accumulates mutations. 
The number of possible states of the tumor is thus $2^d$.
In typical applications one is interested in probability distributions over this state space that are far from stationary. Extrapolating the future course of a given tumor requires the computation of transient distributions.  
However, since the age of a tumor and thus the time point of an observation of the Markov chain is generally unknown, we study the time-marginal distribution, which is defined as the solution of a certain linear system.
Today at least $d = 299$ genes are known to drive tumor progression, cf.~\cite{Bailey2018}, i.e., a single distribution would require the storage of more entries than there are atoms in the visible universe.
This phenomenon is called \emph{state-space explosion}~\cite{Buchholz2007} and renders classical methods for calculating or storing distributions infeasible.
This problem is also known in other application areas, e.g., chemical-reaction networks~\cite{Anderson2010}, the chemical master equation~\cite{Kazeev2014}, Hamiltonian dynamics~\cite{Lubich2016}, queuing networks~\cite{Chan1987}, or stochastic neural networks~\cite{Yamanaka1997}.

To overcome the state-space explosion, we describe such Markov chains via \emph{Stochastic Automata Networks}~\cite{Plateau2000} which provide a sparse representation of the infinitesimal generator as a short sum of Kronecker products.
The main idea is to split up a large Markov chain into smaller ones, the so-called \emph{automata}, which are interconnected.
In tumor progression, for example, each automaton represents a driver mutation and encodes whether it is present or absent.
When describing transitions, there are effects between individual automata, e.g., some mutations may favor or inhibit the occurrence of others.
This leads to entries within the sparse representation of the generator that are functions of the current state.
The explicit evaluation of these functions at each state would result in exponential cost.
However, if these functions are separable, the structure as a short sum of Kronecker products without functional entries can be preserved and the exponential cost avoided.
Assuming this separability, we represent the generator and also the operator of the system as a short sum of Kronecker products.
We make the common assumption that the Markov chain starts always at a certain state, which implies that the right-hand side of the system is a Kronecker product.

To determine the marginal distribution, i.e., to solve the linear system, in a sparse way we exploit this structure and turn to so-called \emph{low-rank tensor formats}.
These are known to handle the state-space explosion by reducing the exponential cost in the number of automata to linear cost.
In the context of low-rank tensors, a sum of Kronecker products is referred to as \emph{CANDECOMP/PARAFAC~(CP) format}~\cite{Carroll1970,Harshman1970}.
For large Markov chains, the CP format has been used extensively to approximate, in particular, stationary distributions, e.g., in~\cite{Kulkarni2011}, and to derive conditions for their existence, e.g., in~\cite{Fourneau2008}.
However, the problem of finding a CP decomposition of a given tensor is NP-hard in general, cf.~\cite{Hastad1990}.
Therefore we focus on tree tensor networks, especially on the \emph{tensor-train}~\cite{Oseledets2009,Ostlund1995,White1992,Loan2008} and the\emph{ hierarchical Tucker format}~\cite{Grasedyck2010,Hackbusch2009}.

There are several known methods to determine an approximate solution of a linear system using low-rank tensor formats, see, e.g.,~\cite{Grasedyck2013} for an overview.
These methods can roughly be divided into optimization-based approaches and iterative procedures including truncation.
Truncation, i.e., approximation of a tensor with one of lower rank, is needed to keep computational and storage costs low during the iteration.
The tensor-train format was successfully used for the computation of, e.g., transient distributions~\cite{Johnson2010}, mean times to failure~\cite{Masetti2019}, and stationary distributions~\cite{Kressner2016,Kressner2014}.
There mainly optimization-based approaches were presented.
In~\cite{Buchholz2016} the hierarchical Tucker format was applied to reduce the storage cost for distributions and the computational cost for performing basic operations.
Continuing in~\cite{Buchholz2017} adaptive truncation strategies for the computation of stationary distributions using iterative methods were presented.

To allow for a probabilistic interpretation of an approximation, the latter should be a probability distribution itself, i.e., all entries should be non-negative, and the sum of all entries should be equal to one.
Similar to~\cite{Kressner2016,Buchholz2017,Kressner2014,Kulkarni2011,Masetti2019} we neglect the non-negativity to ensure an error-controlled approach, see, e.g.,~\cite{Kim2014} for an overview.
We focus on an iterative solver to determine time-marginal distributions which guarantees that the sum of all entries of the approximation is one.
In~\cite{Kressner2014} a power iteration based on a formulation of the stationary solution as an eigenvalue problem was derived, where after each application of the operator in the power iteration the current approximation was rescaled to ensure that it sums up to one.
However, the marginal distribution we are interested in is not stationary, i.e., the right-hand side of the linear system is non-zero.
Rescaling the current approximation is not an option in this case since the final result would then solve the rescaled system rather than approximate the marginal distribution.
Alternatively, one could formulate the system as an eigenvalue problem, but this is not straightforward because it is unclear how the right-hand side should be modified for this purpose.
Other popular iteration methods, e.g., projection methods, which are compatible with low-rank tensors do not allow us to implement our normalization condition either since the convergence of these methods cannot be guaranteed if rescaling takes place after every step.

To overcome this problem, we derive a novel iterative method based on the Neumann series~\cite{Dubois1979} and the uniformization method~\cite{Grassmann1977} using low-rank tensor formats. We verify that this method results in an approximation that sums up to one and prove its convergence.
In our numerical experiments, we focus on the concept of \emph{Mutual Hazard Networks}~\cite{Schill2019} for tumor progression.
Our experiments illustrate that the marginal distribution can be approximated by low-rank tensors using the new algorithm.

This work is organized as follows.
In~\cref{section:problem_statement} we derive the linear system that defines the time-marginal distribution of a continuous-time Markov chain. We represent its operator and right-hand side in a sparse way based on Stochastic Automata Networks.
In~\cref{section:low_rank_method} we introduce the concept of tensors and tree tensor networks. We present the tensor-train and the hierarchical Tucker format.
Using these formats we derive an iterative method and prove its convergence.
In~\cref{section:numerical_experiments} we perform numerical experiments based on the concept of Mutual Hazard Networks.
In~\cref{section:conclusion} we conclude.

\section{Statement of the problem}
\label{section:problem_statement}

First we explain the main problem we address.

\subsection{Time-marginal distribution}
\label{subsection:time_marginal_distribution}

A continuous-time Markov chain is defined by its state space $\Space$, its infinitesimal generator $\Q$ and an initial distribution $\pN$. 
In this paper we assume that the state space $\Space$ is discrete.
The generator is an operator $\Q \in \R^{\Space \times \Space}$ which stores the rates of transition from state $x$ to another state $y$ in $\Q\left[y, x\right]$ and the rates of staying in a state $x$ in $\Q\left[x, x\right]$.
The probability distribution at time $t \geq 0$ is defined by 
\begin{align}
	\label{eq:transient_distribution_with_time}
		\frac{\mathrm{d}}{\mathrm{d} t}\p(t) = \Q  \pN \qquad \text{ or } \qquad \p(t) = \exp(t \Q) \pN .
\end{align}
We make the common assumption that every trajectory starts at the same state, i.e., the initial distribution $\pN \in \R^{\Space}$ is a canonical unit vector.

Here, we assume that the observation time $t$ is unavailable, such as, e.g., in tumor progression modeling, and thus must be treated as a random variable.
Therefore, we are interested in a so-called \emph{time-marginal distribution} $\p$ which is independent of the time $t$ and which we will call marginal distribution for brevity.
Each entry $\p[x]$ of the marginal distribution indicates the probability of observing a state $x \in \Space$ at a random time.
We follow the common assumption that the sampling time is an exponentially distributed random variable with mean $1$, i.e., $t \sim \operatorname{Exp}[1]$. Similar approaches can be found, e.g., in~\cite{Beerenwinkel2009,Hjelm2006,Schill2019}.
Thus, we have 
\begin{align}
	\label{eq:marginal_distribution_integral}
	\p = \int\limits_{0}^{\infty} \exp(-t) \p(t)~ \mathrm{d} t = \int\limits_{0}^{\infty} \exp(t \left(\Q - \Id \right)) \pN~ \mathrm{d} t .
\end{align}
In the following we assume that the spectrum $\sigma\left(\Q - \Id\right)$ is contained in the negative complex half plane $\C^{-} = \{z \in \C ~ \vert ~ \operatorname{Re}\left(z\right) < 0 \}$. In this case the improper integral exists and
\begin{align}
	\label{eq:marginal_distribution}
	\p = \left(\Id - \Q\right)^{-1} \pN.
\end{align}
Hence, the marginal distribution $\p$ is defined as the unique solution of a linear system, since the operator $\Id -  \Q$ is regular.
To overcome the state-space explosion, we focus on the concept of \emph{Stochastic Automata Networks}~\cite{Plateau2000} which offer a sparse representation of the infinitesimal generator $\Q$.

\subsection{Stochastic Automata Networks}
\label{subsection:stochastic_automata_networks}

For a continuous-time Markov chain of interacting processes the discrete state space factorizes in a natural way into the state spaces of the individual processes.
Each process is itself a Markov chain over its own state space $\Space_{\mu}$ and is called a \emph{stochastic automaton} $\A_{\mu}$.
The set $\{\A_1, \dots, \A_d\}$ of $d$ stochastic automata is called a \emph{Stochastic Automata Network}, cf.~\cite{Plateau2000}.
The full state space $\Space$ consists of all possible combinations of states in each automaton, i.e., for $n_{\mu} = \vert \Space_{\mu} \vert$ and $n = \max_{\mu} n_{\mu}$ it is given by
\begin{align}
\label{eq:state_space}
	\Space = \bigtimes\limits_{\mu = 1}^d \Space_{\mu} \qquad \text{with} \qquad \vert \Space \vert = \prod\limits_{\mu = 1}^d n_{\mu} = \mathcal{O}\left(n^d\right).
\end{align}
Each state $x \in S$ specifies its state $x[\mu] \in \Space_{\mu}$ in each automaton $\A_{\mu}$.

We now consider transitions between states in the full state space $\Space$, which are given by one or more transitions between states within the individual state spaces $\Space_{\mu}$.
Such transitions are called \emph{functional}/\emph{constant} and \emph{synchronized}/\emph{local} according to the following definitions.
\begin{itemize}
	\item[(a)] 	A transition is called functional if the transition rate is a function of the current state in $\Space$. Otherwise it is called constant.
	\item[(b)] A transition in a state space $\Space_{\mu}$ may force or prevent transitions in another state space $\Space_{\nu}$ at the same time. Such transitions in $\Space$ are called synchronized. Otherwise they are called local.
\end{itemize}
It is known that the infinitesimal generator $\Q$ of a Stochastic Automata Network with $d$ automata can be represented as a sum of $d$ Kronecker products, cf.~\cite{Plateau2000}.
To this end,  the infinitesimal generator is split into two parts describing local and synchronizing transitions, respectively, i.e., $\Q = \Q_{\text{local}} + \Q_{\text{sync}} $.
Since the local transitions only change the state in their corresponding automaton, we represent the local part as
\begin{align}
\label{eq:local_transition_part}
	\Q_{\text{local}}  = \sum\limits_{\nu = 1}^d ~ \bigotimes\limits_{\mu < \nu} \Id_{\Space_\mu} \otimes ~\Q_{\text{local}} ^{(\nu)} \otimes \bigotimes\limits_{\mu > \nu} \Id_{\Space_\mu}
\end{align}
where $\Id_{\Space_{\mu}}$ denotes the identity in $\Space_{\mu}$ and $\Q_{\text{local}} ^{(\nu)}$ describes the local transition in $\Space_{\nu}$.
The synchronized transitions are represented similarly. Instead of one term of Kronecker products, two terms are needed for each synchronized transition, cf.~\cite{Plateau2000}.
The whole infinitesimal generator for $d$ automata and $s$ synchronized transitions is given by
\begin{align}
\label{eq:general_cp_transiton_matrix}
	\Q = \Q_{\text{local}} + \Q_{\text{sync}} = \sum\limits_{\nu = 1}^{2 s + d} ~ \bigotimes\limits_{\mu = 1}^d \Q_{\nu}^{(\mu)}.
\end{align}
Instead of $\vert \Space \vert \cdot \vert \Space \vert = \mathcal{O}\left(n^{2d}\right)$ entries, one only has to store $\left(2 s + d\right) d n^2$ entries.
Following the current state of the art in tumor progression models we only consider local transitions from now on.

\subsection{Description of functional transitions}
\label{subsection:parameterization}

For Markov chains describing interacting processes, i.e., automata, we are interested in  functional transitions, for which the entries of $\Q_{\nu}^{(\mu)}$ are functions.
Every time such an entry is needed one would have to evaluate the function, which requires knowledge of  the states of all automata.

To analyze this problem more precisely, we concentrate on Markov chains satisfying the following assumptions:
\begin{enumerate}
	\item \label{ass:(i)}
	The continuous-time Markov chain can be represented as a Stochastic Automata Network $\{\A_1, \dots, \A_d\}$ with $d$ automata.
	Each automaton $\A_{\mu}$ has a state space $\Space_{\mu}$ with $\vert \Space_{\mu} \vert = n_{\mu}$, and the state space of the network is given by $S = \bigtimes_{\mu = 1}^d \Space_{\mu}$.
	\item \label{ass:(ii)}
	There are only local and functional transitions.
	\item \label{ass:(iii)}
	For each automaton $\A_{\mu}$ there exists an ordering of the states $i_1 < i_2 < \cdots < i_{n_{\mu}}$ such that there is no transition from state $i$ to state $j$ for all $i > j$. We denote the set of possible transitions in $\Space_{\mu}$ by $T_{\mu} = \{ \left(i, j\right) ~ \vert ~ i < j \in \Space_{\mu} \}$.
	\item \label{ass:(iv)}
	The transition rate from state $x$ to $y$ in $\Space$, where $x$ and $y$ differ only in automaton $\A_{\nu}$ and satisfy $x[\nu] < y[\nu]$, is separable and can be represented with parameters $ \PM_{ \left(x[\nu], y[\nu] \right), x[\mu]} \geq 0$ by
	\begin{align*}
		Q[y, x] = 	\prod\limits_{\mu = 1}^d \PM_{ \left(x[\nu], y[\nu] \right), x[\mu]} .
	\end{align*}
	\item \label{ass:(v)}
	The time $t$ is an exponentially distributed random variable with mean $1$, i.e., $t \sim \operatorname{Exp}[1]$. 
\end{enumerate}
In the following we discuss these assumptions in more detail.
Due to assumption~\ref{ass:(ii)} and~\ref{ass:(iii)}, all possible transitions are covered in assumption~\ref{ass:(iv)}. 
All other transition rates are zero.
Each diagonal entry of the generator is given by minus its off-diagonal column sum, since each column sum is zero.
Hence, the whole infinitesimal generator is already defined by the parameters $\PM_{\left(i_{\nu}, j_{\nu}\right), i_{\mu}} \geq 0$ for all states $i_{\mu} \in \Space_{\mu}$ and $\left(i_{\nu},  j_{\nu}\right) \in T_{\nu}$, and we denote the corresponding generator as $\QPM$.
The parameter $\PM_{\left(i_{\nu}, j_{\nu}\right), i_{\mu}}$ with $\nu \neq \mu$ can be interpreted as a direct multiplicative effect of the state $i_{\mu} \in \Space_{\mu}$ on the local transition from $i_{\nu}$ to $j_{\nu}$ in $\Space_{\nu}$, while $\PM_{\left(i_{\nu}, j_{\nu}\right), i_{\nu}}$ is a  baseline rate of transition from $i_{\nu}$ to $j_{\nu}$ in $\Space_{\nu}$.
If a state $i_{\mu} \in \Space_{\mu}$ has no direct effect on the transition from $i_{\nu}$ to $j_{\nu}$ in $\Space_{\nu}$, then $\PM_{\left(i_{\nu}, j_{\nu}\right), i_{\mu}} = 1$.
As usually assumed, an automaton does not affect many others directly, and therefore many multiplicative effects, i.e., parameters $\PM_{\left(i_{\nu}, j_{\nu}\right), i_{\mu}}$, are equal to one.
This corresponds to a sparsity of the interactions between the automata.
Assumption~\ref{ass:(iv)} allows us to separate the functional transition rates and to transfer the parameters $\PM_{\left(i_{\nu}, j_{\nu}\right), i_{\mu}}$ to the matrices corresponding to automaton $\A_{\mu}$ within each Kronecker product in~\cref{eq:local_transition_part}.
Thus, we represent all functional transitions within the Kronecker product structure without any functional entries:
\begin{align}
\label{eq:Q_CP}
	\QPM &= \sum\limits_{\nu = 1}^d ~ \sum\limits_{\left(i,j\right) \in T_{\nu}} ~ \bigotimes\limits_{\mu = 1} ^d \QCore
	\\
	\text{ with } 	 \qquad
	\QCore &= 
		\begin{cases}
			\diag\left(\PM_{\left(i,j\right), \eta} ~\vert~ \eta \in \Space_{\mu} \right) & \text{if } \mu \neq \nu,
			\\[0.2cm]
			\PM_{\left(i,j\right), i} \left(\Delta_{j,i}^{\nu} - \Delta_{i,i}^{\nu}\right) & \text{if } \mu = \nu
		\end{cases}
	\nonumber
\end{align}
where $\Delta_{j,i}^{\nu} \in \R^{\Space_{\nu} \times \Space_{\nu}}$ has only one non-zero entry $\Delta_{j,i}^{\nu}[j,i] = 1$.
By assumption~\ref{ass:(iii)}, the generator $\QPM$ is similar to a lower triangular matrix where the similarity transformation is obtained by permuting states following the order of each automaton. 
Hence, the spectrum consists of all diagonal entries
\begin{align*}
	\sigma\left(\QPM - \Id\right) = \Bigg\{-1 - \sum\limits_{\nu = 1}^d ~ \sum\limits_{\left(x[\nu], j\right) \in T_{\nu} } ~ \prod\limits_{\mu = 1}^d 	\PM_{\left(x[\nu], j\right), x[\mu]} ~  \Bigg\vert ~x \in S\Bigg\} \subseteq \R_{\leq - 1}
\end{align*}
and the improper integral in~\cref{eq:marginal_distribution} exists.

\subsection{Mutual Hazard Networks}
\label{subsection:mutual_hazard_networks}

A class of Markov chains satisfying assumptions~\ref{ass:(i)} to~\ref{ass:(v)} occurs in the context of tumor progression, the \emph{Mutual Hazard Network} model~\cite{Schill2019}.
There tumor progression is modeled as a Stochastic Automata Network over a discrete state space $\Space$, and a state $x \in \Space$ represents the genotype of a tumor.
Each automaton $\A_{\mu}$ represents a genomic event such as a point mutation, a copy number alteration, or a change in DNA methylation.
The state space of each automaton $\A_{\mu}$ is given by $ \Space_{\mu} = \{0,1\}$, where $x[\mu] = 0$ indicates that the genomic event has not occurred yet and $x[\mu] = 1$ that it has.
Hence, the state space of the Markov chain is given by $\Space = \bigtimes_{\mu = 1}^d \{0,1\}$.
All events are assumed to occur one after another and irreversibly, i.e., there are only local transitions from state $x[\mu] = 0$ to $x[\mu] = 1$ in $\Space_{\mu}$.
If an event has not yet occurred, this has no direct effect on all other events, i.e., $\PM_{\left(0_{\nu}, 1_{\nu}\right), 0_{\mu}} = 1$ for all $\mu \neq \nu$.
A parameter $\PM_{ \left(0_{\nu}, 1_{\nu}\right), 1_{\mu}} > 1$ is interpreted as a promoting effect and $\PM_{\left(0_{\nu}, 1_{\nu}\right), 1_{\mu}} < 1$ as an inhibiting effect.
Furthermore, time is considered as a random variable $t \sim \operatorname{Exp}[1]$, cf.~\cite{Schill2019}.
In summary, all assumptions~\ref{ass:(i)} to~\ref{ass:(v)} are satisfied.
Similar conclusions can be obtained for several types of progression models, e.g.,~\cite{Beerenwinkel2009,Hjelm2006}.

The Mutual Hazard Network with $d$ genomic events corresponds to parameters 
$\PM_{ \left(0_{\nu}, 1_{\nu}\right), 1_{\mu}} > 0$ and $\PM_{\left(0_{\nu}, 1_{\nu}\right), 0_{\mu}} = 1$ for all $\mu, \nu \in \setd$.
Thus, the parameters can be collected in a matrix $\PM \in \R^{d \times d}$ with $\PM[\nu, \mu] = \PM_{\left(0_{\nu}, 1_{\nu}\right), 0_{\mu}}$.
To allow for visualization and interpretation of the network, one assumes sparsity of the Mutual Hazard Network, i.e., many automata or events do not interact directly.
In~\cite{Schill2019} optimal parameters $\PM$ are found using maximum-likelihood estimation for a data set of tumors.
Since the age of a detected tumor is unknown, the likelihood is given by the marginal distribution $\pPM$, i.e., each entry $\pPM [x]$ is the probability under the model that a tumor has the genotype $x$ when it is observed at a random time.
This requires $\pPM$ for fixed parameters in each optimization step.
Due to the state-space explosion, the computation of the marginal distribution in~\cite{Schill2019} based on classical methods was limited to about $d < 25$ automata.

\subsection{Structure of the linear system}
\label{subsection:structure_linear_systems}

As the improper integral in~\cref{eq:marginal_distribution_integral} exists, the time-marginal distribution $\pPM$ depending on parameters $\PM$ is uniquely defined by
\begin{align}
	\label{eq:linear_system}
	\left(\Id - \QPM\right) \pPM = \pN .
\end{align}
Similar to the representation in~\cref{eq:Q_CP}, the identity $\Id  = \bigotimes_{\mu =1}^d \Id_{\Space_{\mu}} \in \R^{\Space \times \Space}$ is a Kronecker product of $d$ smaller identities $\Id_{\Space_{\mu}} \in \R^{\Space_{\mu} \times \Space_{\mu}}$ corresponding to the state spaces $\Space_{\mu}$. The initial distribution $\pN$ is a canonical unit vector and thus a Kronecker product of $d$ canonical unit vectors corresponding to each $\Space_{\mu}$.
Hence, the operator and the right-hand side of~\cref{eq:linear_system} have a representation as a short sum of $d$ Kronecker products, which allows for efficient storage.

\section{Low-rank method}
\label{section:low_rank_method}

We now compute the marginal distribution $\pPM$ in a sparse way.
To avoid losing sparse structures when performing arithmetic operations we keep the Kronecker products unexecuted and understand our operators and distributions as \emph{tensors}.

\subsection{Low-rank tensor formats}
\label{subsection:low_rank_tensors}

We view tensors as multidimensional generalizations of vectors and matrices, i.e., one-dimensional and two-dimensional tensors.
\begin{definition}[tensor]
\label{def:tensor}
	Let $d \in \N$ and $\mathcal{I} = \times_{\mu=1}^d \mathcal{I}_{\mu}$ be a Cartesian product of discrete index sets $\mathcal{I}_{\mu}$.
	An object $\mathcal{B} \in \R^{\mathcal{I}}$ is called a \emph{tensor} of \emph{dimension} $d$.
	Each direction $\mu \in \{1, \dots, d\}$ is called a \emph{mode} of $\mathcal{B}$, and the cardinality of the $\mu$-th index set $\vert \mathcal{I}_{\mu} \vert$ is called the \emph{$\mu$-mode size}.
\end{definition}
In our case, the index set $\mathcal{I}$ corresponds to the state space $\Space = \bigtimes_{\mu = 1}^d \Space_{\mu}$, our distributions $\p$, $\pPM$, $\pN \in \R^\Space$ are tensors of dimension $d$, and the automata $\A_{\mu}$ correspond to the modes with sizes $n_{\mu} = \vert \Space_{\mu} \vert$.

Similar to the reshaping of the distribution vectors into tensors, one can reshape a tensor into a matrix.
This concept, called \emph{matricization}, corresponds to an unfolding of the tensor and is defined, in analogy to~\cite{Grasedyck2010}, as follows.
\begin{definition}[matricization]
	Let $\mathcal{B} \in \mathbb{R}^{\I}$ and $t \subseteq \setd$ with $t \neq \emptyset$ and $s = \setd \setminus t$. 
	The \emph{matricization} of $\mathcal{B}$ corresponding to $t$ is defined as $\mathcal{B}^{(t)} \in \mathbb{R}^{\I_t \times \I_s}$ with $\I_t = \prod\limits_{\mu \in t} \I_{\mu}$ and 
	\begin{align*}
		\mathcal{B}^{(t)} [ (i_{\mu})_{{\mu} \in t},(i_{\mu})_{{\mu} \in s}] = \mathcal{B} [i_1,\dots,i_d]
	\end{align*}
	for all $i=(i_{\mu})_{{\mu} \in \setd} \in \I$. In particular $\mathcal{B}^{(\setd)} \in \mathbb{R}^{\I}$.
\end{definition}
In the language of tensors, the structure of  $\QPM$ in~\cref{eq:Q_CP} can be generalized to the so-called \emph{CANDECOMP/PARAFAC (CP) format} introduced in~\cite{Carroll1970,Harshman1970}.
\begin{definition}[CP format]
\label{def:CP_format}
	A tensor  $\mathcal{B}\in \R^{\mathcal{I}}$ has a \emph{CP representation} if there exist $b_{\nu}^{\left( \mu\right)} \in \R^{\mathcal{I}_{\mu}}$ such that
	\begin{align}
	\label{eq:CP_representation}
		\mathcal{B} = \sum\limits_{\nu = 1}^r \bigotimes\limits_{\mu = 1}^d b_{\nu}^{\left(\mu \right)}.
	\end{align}
	Then $r \in \N_0$ is called the \emph{CP representation rank}, and the  $b_{\nu}^{\left( \mu\right)}$ are called the \emph{CP cores} of $\mathcal{B}$.
\end{definition}
The infinitesimal generator $\QPM$ in~\cref{eq:Q_CP} has CP representation rank
$\sum_{\nu = 1}^d \vert T_{\nu}\vert \leq d n^2$, while the identity $\Id$ as well as the right-hand side $\pN$ have rank $1$.
A core advantage of the CP format is the data sparsity in case of small representation
rank $r$: The representation~\cref{eq:CP_representation} of a tensor $\mathcal{B} \in \R^{\mathcal{I}}$ has storage cost in $\mathcal{O}\bigl( r \sum_{\mu = 1}^d n_{\mu} \bigr) = \mathcal{O}(rdn)$ in contrast to $\mathcal{O}\bigl(\prod_{\mu = 1}^d n_{\mu}\bigr) = \mathcal{O}\left(n^d\right)$ for $n = \max_{\mu} n_{\mu}$.

The problem of finding conditions that guarantee the existence of a low-rank approximation of a given tensor is an active research topic of its own, see, e.g.,~\cite{Bachmayr2017,Kressner2016,Plateau2007}.
Since this goes beyond the scope of this article, we concentrate on numerical experiments which indicate that the marginal distribution $\pPM$ can be approximated with low rank.

To this end, we have to solve a linear system whose operator and right-hand side have a CP representation.
For an operator with CP rank $r > 1$ it is unknown how to calculate its inverse analytically. 
Hence, we need a solver and arithmetic operations to compute the solution numerically.
Most of these arithmetic operations lead to an increase in the representation rank and thus in storage.
To avoid such an increase one would like to \emph{truncate} a given CP tensor to lower rank, i.e., approximate the tensor with one of lower rank.
As the set of tensors with CP rank at least $r$ is not closed for $d > 2$, low-rank approximation within the CP format is an ill-posed problem, cf.~\cite{Silva2008}.
To overcome this drawback, we use the following low-rank tensor formats which allow for truncation in an error-controlled way.

\subsection{Tree tensor network graphs}
\label{subsection:tree_tensor_networks}

\begin{figure}[tbhp]
	\centering
	
	\includegraphics[width = 0.3\textwidth]{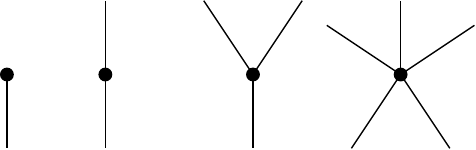}
	
	\caption{\footnotesize{Tensor network graphs: vector, matrix, three-dimensional and five-dimensional tensor (from left to right).}}
	\label{fig:tensor_network_graphs}
\end{figure}

\begin{figure}[tbhp]
	\centering
	
	\includegraphics[width = 0.3\textwidth]{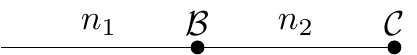}
	
	\caption{\footnotesize{Matrix-vector product $\mathcal{B} \cdot \mathcal{C} \in \R^{\I_1}$ as a tensor network graph with $\mathcal{B} \in \R^{\I_1 \times \I_2}$, $\mathcal{C} \in \R^{\I_2}$, and $n_{\mu}=\vert \I_{\mu} \vert$.}}
	\label{figure:matrix_vector_product}
\end{figure}

\begin{figure}[tbhp]
	\centering

	\includegraphics[width = 0.3\textwidth]{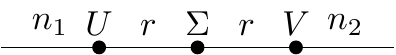}
	
	\caption{\footnotesize{Singular value decomposition $\mathcal{B} = U \Sigma V^T \in \R^{\I_1 \times \I_2}$ as a tensor network graph with matrix rank $r$ and $n_{\mu}= \vert \I_{\mu} \vert$.}}
	\label{figure:SVD}
\end{figure}

Tensors and many interactions between them can be illustrated via an undirected graph as a \emph{tensor network}.
Each vertex of a tensor network graph represents a particular tensor. The dimension of this tensor is given by the number of its edges and half-edges, where a half-edge is an edge that is connected to only one vertex.~\Cref{fig:tensor_network_graphs} shows some examples of tensor network graphs.
An edge between two vertices can be interpreted as a \emph{contraction} over the corresponding mode, i.e., one sums over the index set of this mode.
The resulting network can again be viewed as a tensor whose dimension is given by the number of half-edges.
For example, the contraction of  a two-dimensional tensor $\mathcal{B} \in \R^{\I_1 \times \I_2}$ and a one-dimensional tensor $\mathcal{C} \in \R^{\I_2}$ is the matrix-vector product $\mathcal{B} \cdot \mathcal{C} \in \R^{\I_1}$ with
\begin{align*}
	\left(\mathcal{B} \cdot \mathcal{C} \right)[x] = \sum\limits_{y \in \I_2} \mathcal{B}[x, y] \mathcal{C}[y]
\end{align*}
for all $x \in \I_1$. The corresponding tensor network graph is shown in~\cref{figure:matrix_vector_product} where the mode sizes are denoted by $n_{\mu} =\vert \I_{\mu} \vert$.
This concept can be extended to higher dimensions.
For a two-dimensional tensor $\mathcal{B} \in \R^{\I_1 \times \I_2}$, i.e., a matrix, one can use these networks to illustrate the \emph{singular value decomposition} $\mathcal{B} = U \Sigma V^T$ as shown in~\cref{figure:SVD}. There, the mode sizes $n_{\mu}= \vert \I_{\mu} \vert $ corresponding to the half-edges equal the original ones of $\mathcal{B}$, and the mode sizes corresponding to the edges are given by the matrix rank $r$.\footnote{For a two-dimensional tensor the matrix rank equals the minimal CP representation rank, cf.~\cite{Hackbusch2012}.}
Again this can be extended to higher dimensions, i.e., higher-dimensional tensors can also be represented as contractions of auxiliary tensors.
In analogy to the two-dimensional case, the mode sizes of the half-edges correspond to those of the original tensor, and the mode sizes of the edges generalize the concept of matrix rank.
This generalized rank allows for an error-controlled way of truncation in tree tensor networks.
Several tensor formats based on tree tensor networks are known.
Here, we focus on the \emph{tensor-train format}~\cite{Oseledets2009,Ostlund1995,White1992,Loan2008} and the \emph{hierarchical Tucker  format}~\cite{Grasedyck2010,Hackbusch2009}.

\subsection{Tensor-train and hierarchical Tucker  format}
\label{subsection:TT_and_HT_format}

The tensor-train format was first introduced to the numerical analysis community in~\cite{Oseledets2009}. It is also known in other areas as \emph{matrix product states}~\cite{Ostlund1995,White1992} or as \emph{linear tensor network}~\cite{Loan2008}.

\begin{figure}[tbhp]
	\centering
	
	\includegraphics[width = 0.25\textwidth]{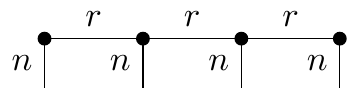}
	
	\caption{\footnotesize{Tensor-train representation of a tensor with dimension $d = 4$,  mode sizes all equal to $n$ and representation rank $\mathbf{r}=\left(r,r,r\right)$.}}
	\label{figure:TT}
\end{figure}

A tensor of dimension $d$ in the tensor-train format is represented as $d-1$ contractions of $d$ so-called \emph{core tensors}. 
Each core tensor represents one mode and is connected to its forward and backward neighbors. The core tensors themselves are of dimension three, or two in case of the first and the last mode.
For a given tensor-train representation, its rank is defined as the tuple of all mode sizes corresponding to edges.~\cref{figure:TT} shows an illustrative example.
Instead of the original tensor one needs to store the core tensors. Thus, the required storage cost for a tensor-train representation with dimension $d$,  mode sizes all equal to $n$ and representation rank component-wise bounded by $r$ is in $\mathcal{O}\left( d n  r^2 \right)$.
In our case the dimension of the distribution tensors is equal to the number $d$ of automata, and the mode sizes are given by the number of states $n_{\mu}$ in each automaton.
Note that in general the representation rank can depend on the ordering of the modes/automata, cf.~\cite{Grasedyck2011,Oseledets2009}.
\\

\begin{figure}[tbhp]
	\centering
	
	\includegraphics[width = 0.25\textwidth]{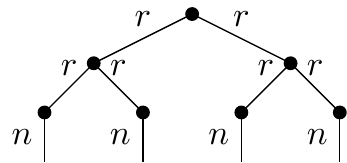}
	
	\caption{\footnotesize{Hierarchical Tucker representation of a tensor with dimension $d = 4$, mode sizes all equal to $n$ and representation rank $\mathbf{r} = \left(r\right)_{t \in \mathcal{T}}$ using a balanced tree $\mathcal{T}$.}}
	\label{figure:HT}
\end{figure}

The hierarchical Tucker format was introduced in~\cite{Hackbusch2009} and further analyzed in~\cite{Grasedyck2010}.
As the name suggests, the hierarchical Tucker format is based on a hierarchical binary tree structure defining the arrangement of the auxiliary tensors. 
Similar to the tensor-train format there are auxiliary tensors each of which represents one mode of the original tensor. These two-dimensional tensors, the so-called \emph{frames}, build the leaves of the binary tree. The two-dimensional tensor at the root of the tree and the three-dimensional tensors at all other levels are called \emph{transfer tensors}. For a given hierarchical Tucker representation, its rank $\mathbf{r}$ is defined as the tuple of all mode sizes corresponding to the edges of the tree, i.e., $\mathbf{r}= \left(r_t\right)_{t \in \mathcal{T}}$ where $\mathcal{T}$ denotes the given tree.~\cref{figure:HT} shows an illustrative example.
Note that in the hierarchical Tucker format we have many options for the structure of the binary tree and that the rank depends on our choice, cf.~\cite{Grasedyck2011}.
Instead of the original tensor the frames and transfer tensors have to be stored. Thus, the storage cost for a hierarchical Tucker representation of a $d$-dimensional tensor with maximum mode size $n$ and a representation rank component-wise bounded by $r$ is in $\mathcal{O}\left( d n r + d  r^3 \right)$.
In our case the dimension of the distribution tensors is again equal to the number $d$ of automata, and the mode sizes are the number of states $n_{\mu}$ for each automaton.
\\
 
Thus, in both the tensor-train and the hierarchical Tucker format, the storage cost grows only linearly in the dimension $d$. This allows for efficient storage and prevents the state-space explosion provided we have low ranks.
It is easily possible to convert a CP representation of rank $r$ to a tensor-train or a hierarchical Tucker representation, where all rank components are bounded by $r$, cf.~\cite{Hackbusch2012}.
Similar to the CP format both formats provide arithmetic operations like addition, application of operators, or scalar products.
Since the storage cost increases with the square or cube of the rank and arithmetic operations also increase the rank, we need a way to truncate a tensor to lower rank in an error-controlled way.
In tree tensor networks we use a generalization of the singular value decomposition to truncate in a quasi-optimal way, see~\cite{Oseledets2009} for the tensor-train and~\cite{Grasedyck2010} for the hierarchical Tucker format.
There it is shown that the error in the Frobenius norm caused by the truncation of a tensor $\mathcal{B} \in \R^{\I}$ to rank $\mathbf{r} = \left(r_{t}\right)_{t \in \mathcal{T}}$ for a given tree $\mathcal{T}$ is bounded by
\begin{align}
	\label{eq:truncation_error}
	\Vert\mathcal{B} - \tau_{\mathbf{r}}(\mathcal{B}) \Vert^2 \leq \sum\limits_{t \in\mathcal{T}} \sum\limits_{m > r_{t}} \sigma_{t, m}^2 \leq~ C d ~\Vert\mathcal{B} - \mathcal{B}^{\operatorname{best}} \Vert^2 
\end{align}
where $\tau_{\mathbf{r}}$ denotes the truncation operator, $\sigma_{t, m}$ is the $m$-th singular value of the matricization $\mathcal{B}^{\left(t\right)}$, $\mathcal{B}^{\operatorname{best}}$ is a best rank-$\mathbf{r}$ approximation, and $C < 2$ is a small constant.
According to~\cref{eq:truncation_error}, a steep decay in the singular values of the corresponding matricization leads to a small relative truncation error, i.e., the tensor can be well approximated with one of low rank.

Inspired by~\cite{Hackbusch2012},~\cref{tab:tensor_arithmetic} summarizes some arithmetic operations and their cost in the tensor-train and the hierarchical Tucker format, see~\cite{Grasedyck2010,Hackbusch2012,Oseledets2009}.~\Cref{tab:tensor_arithmetic} shows that using those low-rank tensor formats reduces the exponential computational cost in the dimension $d$ to linear cost.
\begin{table}[tbhp]
	\captionsetup{justification=centering}
	{\footnotesize
			\caption{\footnotesize \\ Operations and their cost in tensor-train $\left(\TT\right)$ and hierarchical Tucker $\left(\HT\right)$  format for $d$-dimensional tensors with maximum mode size $n$ and rank component-wise bounded by $r$.}
			\label{tab:tensor_arithmetic}
		\begin{center}
			\begin{tabular}{lll}
				\hline
				Operation & Cost in $\TT$ & Cost in $\HT$  
				\\ 
				\hline
				Storage & $\mathcal{O}( d n r^2)$ &$\mathcal{O}( d r^3 + d n r)$ 
				\\
				Addition & $\mathcal{O}( d n r^3)$  & $\mathcal{O}( d n  r^2 + d r^4)$ 
				\\
				Evaluation & $\mathcal{O}( d r ^2)$  &  $\mathcal{O}( d r^3 )$
				\\
				Scalar product& $\mathcal{O}( d n r^3)$& $\mathcal{O}( d n r^2 + d r^4)$ 
				\\
				Operator application & $\mathcal{O}( d n^2 r^4)$ & $\mathcal{O}( d n^2 r)$ 
				\\
				Truncation  &  $\mathcal{O}( d n r^3 )$ & $\mathcal{O}( d n r^2 + d r^4 )$ 
				\\ 
				\hline
			\end{tabular}
	\end{center}
	}
\end{table}
Since the hierarchical Tucker format is based on a binary tree, it is possible to do many of the arithmetic operations and especially the costly truncation in parallel level by level of the tree as described in~\cite{Grasedyck2017}.
Thus, for a balanced tree and dimension $d$ the runtime can often be reduced from $\mathcal{O}\left(d\right)$ to $\mathcal{O}\left(\log(d)\right)$, cf.~\cite{Grasedyck2017}.

\subsection{Low-rank method}
\label{subsection:low_rank_method}

We now make use of these low-rank tensor formats to approximate the marginal distribution $\pPM$ as the solution of~\cref{eq:linear_system}.
Since the exact solution $\pPM$ is a probability distribution, its entries sum up to one, i.e., it fulfills
\begin{align}
\label{eq:probability_condition}
	\langle \e , \pPM \rangle = 1
\end{align}
where $\e \in \R^S$ is the tensor of all ones.
To allow for a probabilistic interpretation of an approximation of $\pPM$, we need to treat~\cref{eq:probability_condition} as an additional constraint.
To this end, we now present an iterative method based on the Neumann series~\cite{Dubois1979} and the uniformization method~\cite{Grassmann1977}.

Let $\gamma \geq \max_{x \in \Space} \vert \QPM[x,x] \vert > 0$ be a bound on the diagonal entries of the infinitesimal generator, then
\begin{align*}
	\pPM = \left(\Id - \QPM\right)^{-1} \pN = \frac{1}{1 + \gamma} \left(\Id - \left(\frac{\gamma}{1 + \gamma} \Id + \frac{1}{1 + \gamma} \QPM\right)\right)^{-1} \pN.
\end{align*}
Owing to the ordering assumption~\ref{ass:(iii)} and the fact that all diagonal entries of $\QPM$ are non-positive, the spectral radius is smaller than one, i.e., 
$\rho\left(\frac{\gamma}{1 + \gamma} \Id + \frac{1}{1 + \gamma} \QPM\right) = \frac{\gamma}{1 + \gamma} < 1$, and so the Neumann series converges.
Hence, we can represent the marginal distribution as
\begin{align}
\label{eq:Neumann_representation}
	\pPM = \frac{1}{1 + \gamma} ~ \sum\limits_{m = 0}^{\infty} ~ \left(\frac{\gamma}{1 + \gamma}\right)^m \left(\Id + \frac{1}{\gamma} \QPM\right)^m \pN.
\end{align}
Alternatively we can derive~\cref{eq:Neumann_representation} using the uniformization method.
Here, the idea is to describe a continuous-time Markov chain by a discrete-time Markov chain with a time increment that is exponentially distributed.
With this interpretation, we write the time-dependent probability distribution in~\cref{eq:transient_distribution_with_time} as
\begin{align*}
	\pPM(t) = \sum\limits_{m = 0}^{\infty} \frac{\left(\gamma t\right)^m}{m!} \exp(-\gamma t) \left(\Pgamma\right)^m \pN \qquad \text{with} \qquad \Pgamma = \Id + \frac{1}{\gamma} \QPM
\end{align*}
for a time $t \geq 0$. 
Note that $\Pgamma $ is the transition probability matrix of a discrete-time Markov chain, since $\gamma$ is a bound on the diagonal entries of $\QPM$.
Marginalization of time similar to~\cref{eq:marginal_distribution_integral} and substitution leads to~\cref{eq:Neumann_representation}:
	\begin{align*}
		\pPM 
		&= \int\limits_{0}^{\infty} \exp(-t) ~ \pPM(t) ~\mathrm{d}t 
		~=~ \sum\limits_{m = 0}^{\infty} ~ \frac{\gamma^m}{m!} ~\int\limits_{0}^{\infty} t^m ~ \exp(- (\gamma + 1 ) t) ~\mathrm{d}t ~ \left(\Pgamma\right)^m \pN
		\\
		&= \sum\limits_{m = 0}^{\infty} ~ \frac{\gamma^m ~ \Gamma(m + 1 )}{ m! \left(1 + \gamma\right)^{m+1}} ~  \left(\Pgamma \right)^m \pN
		~=~ \frac{1}{1 + \gamma} ~ \sum\limits_{m = 0}^{\infty} ~ \left(\frac{\gamma}{1 + \gamma} ~\Pgamma \right)^m \pN
	\end{align*}
where $\Gamma$ denotes the gamma function. A natural approximation based on~\cref{eq:Neumann_representation} would be
\begin{align*}
	\pktilde= \frac{1}{1 + \gamma} ~ \sum\limits_{m = 0}^{k} ~ \left(\frac{\gamma}{1 + \gamma} ~\Pgamma \right)^m \pN
\end{align*}
for $k \in \N$. Based on the properties of the Neumann series this sequence converges linearly to $\pPM$, but $\tilde{\mathbf{p}}^{(k)}$ does not satisfy the normalization condition~\cref{eq:probability_condition} for any $k \in \N$.
As $\Pgamma $ is a transition probability matrix, its application to a probability distribution leads to a probability distribution again satisfying~\cref{eq:probability_condition}, and thus
\begin{align*}
	\langle \e, \sum\limits_{m = 0}^{k} ~ \left(\frac{\gamma}{1 + \gamma} ~ \Pgamma\right)^m \pN \rangle=  \sum\limits_{m = 0}^{k} ~ \left(\frac{\gamma}{1 + \gamma}\right)^m = \frac{\left(1 + \gamma\right)^{k + 1} - \gamma^{k+1}}{\left(1 + \gamma\right)^k}
\end{align*}
for all $k \in \N$.
By scaling each element of the sequence, we obtain
\begin{align}
\label{eq:pk}
	\pk =  \frac{\left(1 + \gamma\right)^k}{\left(1 + \gamma\right)^{k + 1} - \gamma^{k+1}} ~ \sum\limits_{m = 0}^{k} ~ \left(\frac{\gamma}{1 + \gamma}~ \Pgamma\right)^m \pN
\end{align}
for $k \in \N$, which is now an approximation to $\pPM$ that satisfies~\cref{eq:probability_condition}. We prove its linear convergence to $\pPM$ in the following theorem.
\begin{theorem}
\label{lemma:convergence_pk}
	Let $\PM$ be some given parameters, $\pPM$ be the solution of the corresponding linear system~\cref{eq:linear_system}, and $\pk$ be defined by~\cref{eq:pk} for all $k \in \N$. Then $\pk$ converges linearly to $\pPM$ as $k$ approaches infinity, i.e.,
	\begin{align*}
		\lim\limits_{k \to \infty} \pk = \pPM \qquad \text{with} \qquad \big\Vert \pk - \pPM \big\Vert \leq c \cdot \biggl(\frac{\gamma}{1 + \gamma}\biggr)^k \quad \text{for all } k \in \N
	\end{align*}
	where $c =\big \Vert \frac{1}{1 + \gamma} \pN - \pPM \big\Vert + \gamma \Vert \pPM \Vert$.
\end{theorem}
\begin{proof}
	For $\alpha_k = \frac{(1 + \gamma)^k}{(1 + \gamma)^{k + 1} - \gamma^{k+1}} $, $\alpha = \frac{1}{1+\gamma}$ and any $k \in \N_0$ we have
	\begin{align*}
		\vert \alpha_k - \alpha \vert &= \bigg\vert \frac{ \gamma^{k+1}}{(1+\gamma) ((1+\gamma)^{k+1} - \gamma^{k+1})}  \bigg\vert < \frac{\gamma}{1+\gamma} ~ \vert \alpha_{k-1} - \alpha \vert
		\\
		&< \biggl(\frac{\gamma}{1 + \gamma}\biggr)^k ~ \vert \alpha_0 - \alpha \vert = \biggl(\frac{\gamma}{1 + \gamma}\biggr)^{k+1}.
	\end{align*}
	Furthermore we obtain
	\begin{align*}
		\big\Vert \pk - \pPM \big\Vert &\leq \big\Vert \pk - \pktilde \big\Vert + \big\Vert \pktilde - \pPM \big\Vert 
		\\
		&\leq \vert \alpha_k - \alpha \vert \bigg\Vert \sum\limits_{m = 0}^{k} ~ \biggl(\frac{\gamma}{1 + \gamma}~\Pgamma\biggr)^m \pN \bigg\Vert + \biggl( \frac{\gamma}{1+\gamma}\biggl)^k \Vert \alpha \pN - \pPM \Vert 
		\\
		&\leq \bigl( \Vert \alpha \pN - \pPM\Vert + \gamma \Vert \pPM \Vert \bigr) \cdot \biggl(\frac{\gamma}{1 + \gamma}\biggr)^k.
	\end{align*}
	Since $\frac{\gamma}{1 + \gamma} < 1$, the linear convergence follows.
\end{proof}

In order to employ low-rank tensor formats, we extend the iterative method corresponding to~\cref{eq:pk} by truncation.
As shown in~\cite{Hackbusch2008}, a convergent iterative method combined with truncation still converges if the truncation error is sufficiently small. 
However, the truncation could violate condition~\cref{eq:probability_condition}.
Therefore, we have to extend our method to ensure that
\begin{align*}
	\langle \e, \tau\left( \left(\Pgamma\right)^m \pN \right) \rangle = 1 \qquad \text{and} \qquad \langle \e, \tau( \pk )\rangle = 1
\end{align*}
for all $m, k \in \N$ where $\tau$ denotes the truncation operator. We implement these conditions by rescaling after truncation as shown in~\cref{alg:uniformtization}, see lines $7$ and $11$.
The procedure should be continued until the norm of the relative residual is smaller than a given tolerance $\tol >0$.

\begin{algorithm}[tbhp]
	\caption{Low-rank uniformization$\left(\PM, \gamma, \tol \right)$}
	\label{alg:uniformtization}
	\begin{algorithmic}[1]
		\STATE $\Pgamma= \Id + \frac{1}{\gamma} \QPM$
		\STATE $k = 0$
		\STATE $s = 1, \csum = 1$
		\STATE $\pk = \pN, \psum = \pN$
		\WHILE{ $\frac{\Vert \left(\Id - \QPM\right) \pk / \csum -  \pN \Vert}{\Vert \pN\Vert } \geq \tol$ }
		\STATE{ $ \psum = \tau\left(\Pgamma \psum\right)$ }
		\STATE{ $ \psum = \frac{ \psum }{ \langle e, \psum \rangle}$ }
		\STATE{$s = \frac{ \gamma}{ \left( 1 + \gamma\right)} \cdot s$}
		\STATE{$\csum = \csum + s$}
		\STATE{$\mathbf{p}^{(k+1)} = \tau\left(\pk + s~ \psum\right)$}
		\STATE{$\mathbf{p}^{(k+1)} =  \frac{\csum }{ \langle e, \pk \rangle} \cdot \mathbf{p}^{(k+1)}$}
		\STATE{$k = k + 1$}
		\ENDWHILE
		\STATE{$\pk = \frac{\pk}{ \csum}$}
		\RETURN $\pk$
	\end{algorithmic} 
\end{algorithm}
In~\cref{alg:uniformtization} we have to choose a bound $\gamma$ on the diagonal entries of $\QPM$.
The cost for computing all diagonal entries to find the maximum is in $\mathcal{O}(n^d)$.
In case of parameters $\PM$ that vary little,
\begin{align}
\label{eq:upper_bound}
	\gamma  = \sum\limits_{\nu = 1}^d ~ \max\limits_{i_{\nu}} \left( \sum\limits_{\left(i_{\nu}, j_{\nu}\right) \in T_{\nu}} ~ \prod\limits_{\mu = 1}^d ~ \max\limits_{i_{\mu}} ~ \PM_{\left(i_{\nu}, j_{\nu}\right), i_{\mu}} \right)
\end{align}
can be used as an inexpensive bound.
For the computation of~\cref{eq:upper_bound} we need $\mathcal{O}\left(d^2 n^2 T\right)$ comparisons and evaluations\footnote{Strictly speaking, the computational cost of the method is therefore quadratic, rather than linear, in $d$. However, $\gamma$ only needs to be precomputed once, the computational cost of which is negligible.} of the parameter $\PM$ where $T = \max_{\nu} \vert T_{\nu} \vert$ is the maximum number of possible transitions in an automaton.
In the case of strongly varying parameters $\PM$,~\cref{eq:upper_bound} may be a gross overestimation for the diagonal entries of $\QPM$.
The question of how to determine a tighter bound with polynomial effort using low-rank tensor formats will be dealt with in future work.

\section{Numerical experiments}
\label{section:numerical_experiments}

We illustrate our method for the computation of time-marginal distributions in numerical experiments using the following setting based on the model of Mutual Hazard Networks.
\subsection{Setting}
\label{subsection:settings}

We consider $d$ automata and sparse parameters $\PM \in \R^{d \times d}$ with a particular block-diagonal form.
Each block of size $b \times b$ characterizes a subset of $b$ automata which directly affect one another.
There are no direct effects between automata from different blocks. 
In each block, we would like the effects between neighboring automata to be stronger than those between automata that are farther apart.
Therefore, we draw the entries $B[i, j]$ of the blocks from a normal distribution with mean $1$ and standard deviation $\sigma = 2^{- 1 - \vert i - j\vert}$ (and restrict them to $\R^{>0}$).
Two possible examples for $d = 4$ automata and block size $b \in \{2, 4\}$ are given below:
\begin{align*}
		\begin{pmatrix}
		{\color{blue}1.2688 }& {\color{red}1.4585 } & {\color{lightgray}1} & {\color{lightgray}1}  \\
		{\color{red}0.43529  }& {\color{blue}1.4311  }& {\color{lightgray}1}  & {\color{lightgray}1}  \\
		{\color{lightgray}1}  & {\color{lightgray}1}  &{\color{blue}1.1594}  & {\color{red} 0.67308}\\
		{\color{lightgray}1}  & {\color{lightgray}1}  & {\color{red} 0.8916} & {\color{blue}1.1713 }
		\end{pmatrix}
	\quad
		\begin{pmatrix}
			{\color{blue}1.2688 }& {\color{red}1.4585 } & {\color{gray}0.71764} & {\color{lightgray}1.0539}  \\
			{\color{red}0.43529  }& {\color{blue}1.4311  }& {\color{red}0.8916}  & {\color{gray}1.0428}  \\
			{\color{gray}1.4473 }  & {\color{red}1.6924 }  &{\color{blue}1.1594}  & {\color{red} 0.67308}\\
			{\color{lightgray}1.0453}  & {\color{gray}0.99212}  & {\color{red} 0.8916} & {\color{blue}1.1713 }
		\end{pmatrix}
\end{align*}
Parameters of the same color follow the same distribution.

For the application of~\cref{alg:uniformtization} we always choose the bound $\gamma$ as in~\cref{eq:upper_bound}, which for Mutual Hazard Networks can be reduced to
\begin{align}
\label{eq:MHN_upper_bound}
	\gamma =  \sum\limits_{\nu = 1}^d ~ \prod\limits_{\mu = 1}^d \max\{1, \PM[\nu, \mu]\}.
\end{align}
We use a canonical balanced tree tensor network for the application of the hierarchical Tucker format, i.e., the automata are assigned to the leaves following their ordering, see~\cref{figure:sv_d=8_b=4}. We perform all experiments for $100$ randomly generated sample parameters for each combination of block size $b$ and number $d$ of automata.
Unless stated otherwise, we compute low-rank approximations of the marginal distribution $\pPM$ using~\cref{alg:uniformtization} with a maximum relative truncation error of $\eps= 10^{-8}$. The algorithm stops when the norm of the relative residual  is smaller than a tolerance value of $\tol = 10^{-4}$.
The mean values we compute are arithmetic means.

\subsection{Study of singular values}
\label{subsection:singular_values}

One fundamental assumption for solving linear systems using low-rank tensor methods is that a solution can be well approximated by a tensor of low rank. 
According to the error bound in~\cref{eq:truncation_error}, the truncation error is determined by the singular values of the matricization corresponding to the chosen tree network. A fast decay of those singular values indicates that a tensor can be approximated accurately by one of low rank.

\begin{figure}[tbhp]
	\centering
	
	\includegraphics[width = \textwidth]{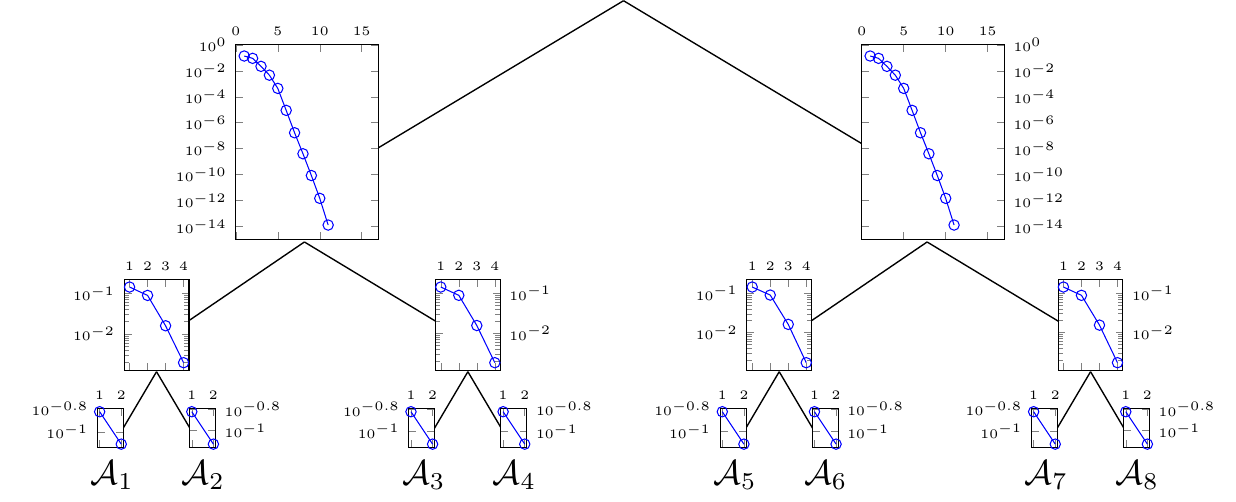}
	
	\caption{\footnotesize{Mean of singular values of matricizations of $\pPM$ for the canonical balanced tree and $100$ sample parameters with $d=8$ automata and block size $b=4$.}}
	\label{figure:sv_d=8_b=4}
\end{figure}

To analyze this issue, we solve~\cref{eq:linear_system} using classical matrix methods of \textsc{MATLAB}~\cite{Matlab2019} and compute the singular values of the corresponding matricization using the \texttt{htucker} toolbox~\cite{Kressner2012}.~\Cref{figure:sv_d=8_b=4} shows the decay of the singular values of each matricization corresponding to the canonical balanced tree for $d = 8$ automata and blocks of size $b = 4$.  For each vertex the semi-logarithmic plot displays the means of the singular values.
We observe that the singular values, especially those close to the root, exhibit an exponential decay.
The two matricizations closest to the root are transposes of each other, and therefore their singular values are identical. The smallest $5$ of the $16$ singular values are indistinguishable from zero and therefore cannot be displayed in the semi-logarithmic plot.
The exponential decay of the singular values indicates that the marginal distribution $\pPM$ is well approximated with low rank.

\subsection{Study of the tree structure}
\label{subsection:study_of_tree}

In general, for tree tensor networks the rank depends on the structure of the tree.
As already observed for large Markov chains using the tensor-train format, a change in the ordering of the automata within the tensor-train tree network can affect the rank, cf.~\cite{Masetti2019}.
In the following, we focus on the hierarchical Tucker format and study how the ordering of the automata in the leaves of the tree affects the low-rank approximability. We preserve the balanced binary structure of the tree because this is advantageous for parallelization, cf.~\cite{Grasedyck2017}.

\begin{figure}[tbhp]
	\centering
	
	\includegraphics[width = \textwidth]{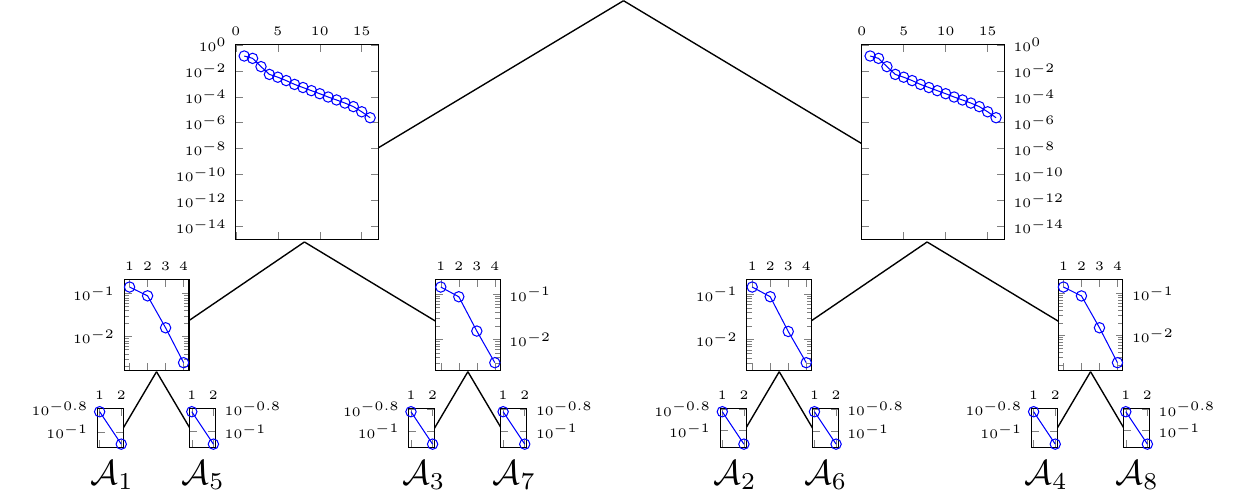}
	
	\caption{\footnotesize{Mean of singular values of matricizations of $\pPM$ for a modified balanced tree and $100$ sample parameters with $d=8$ automata and block size $b=4$.}}
	\label{figure:sv_d=8_b=4_reverse}
\end{figure}

We already studied the singular values of the matricizations corresponding to the canonical balanced tree, see~\cref{figure:sv_d=8_b=4}. We change only the arrangement of the automata in the leaves of the tree and compute the marginal distribution $\pPM$ again.
The decay of the singular values for each matricization is shown in~\cref{figure:sv_d=8_b=4_reverse} where the semi-logarithmic plot at each vertex displays the means of the singular values.
Comparing~\cref{figure:sv_d=8_b=4,figure:sv_d=8_b=4_reverse}, we observe that the choice of the canonical binary tree, i.e., the original ordering in the leaves, results in a significantly faster decay of the singular values close to the root.~\Cref{figure:sv_d=8_b=4_reverse} shows that the singular values closest to the root also have an exponential decrease, but at a much slower rate.
Note that blocks of size $b = 4$ indicate that the automata $\{\A_1, \A_2, \A_3, \A_4\}$ interact directly with one another and are thus highly correlated.
The same holds for the automata  $\{\A_5, \A_6, \A_7, \A_8\}$.
There are no direct interactions between automata of different blocks, i.e., these are only weakly correlated.
Hence, the modified balanced tree separates highly correlated automata, which explains the slower decline of the singular values in level $1$ (level $0$ being the root).
In contrast, the canonical balanced tree separates weakly correlated automata, leading to a faster decline of the singular values.
We will make similar observations when studying the approximation rank.

\subsection{Study of the approximation rank}
\label{subsection:approximation_rank}
As the rank $\mathbf{r} = \left(r_t\right)_{t \in \mathcal{T}}$ in a tree tensor network is a tuple depending on the underlying tree $\mathcal{T}$, we consider the \emph{maximum rank} $\rmax = \max_{t} r_t$ and the \emph{effective rank} $\reff$. The effective rank $\reff$ of a tensor representation is defined such that the storage cost for this representation equals the cost to store one with rank $\mathbf{r} = \left(\reff \right)_{t \in \mathcal{T}}$. Since intuitively a rank should be an integer, $\reff$ is rounded up to the nearest integer.

\begin{figure}[tbhp]
	\begin{minipage}[tbhp]{0.45\textwidth}
		\vspace{0pt}
		\centering

		\includegraphics[width =\textwidth]{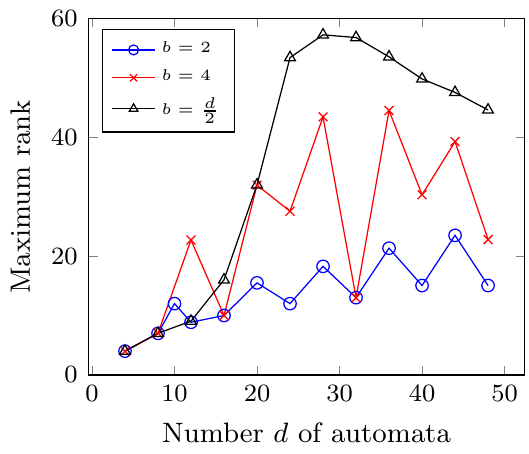}
		
		\subcaption{\footnotesize{Mean of maximum rank.}}
		\label{figure:max_rank_b}
	\end{minipage}
	\hfill
	\begin{minipage}[tbhp]{0.45\textwidth}
		\vspace{0pt}
		\centering

		\includegraphics[width =\textwidth]{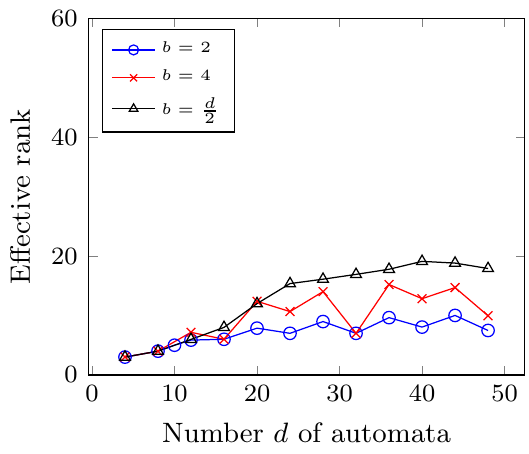}
		
		\subcaption{\footnotesize{Mean of effective rank.}}
		\label{figure:eff_rank_b}
	\end{minipage}
	\caption{\footnotesize{Mean of maximum and effective approximation rank of $\pPM$ as a function of the number $d$ of automata using~\cref{alg:uniformtization} with tolerance $\tol=10^{-4}$ and maximum relative truncation error $\eps = 10^{-8}$ for $100$ sample parameters each.}}
\end{figure}

We compute low-rank approximations of the marginal distribution $\pPM$ using~\cref{alg:uniformtization} as described in~\cref{subsection:settings}.
In the plots we only show the mean values of the ranks, since their variance among the $100$ realizations of $\PM$ was very small.~\Cref{figure:max_rank_b,figure:eff_rank_b} show the maximum and effective approximation rank of $\pPM$ as a function of the number $d$ of automata for blocks of size $b \in \{2, 4, \frac{d}{2}\}$.
We observe that $\pPM$ is approximated to a tolerance of $\tol = 10^{-4}$ with low rank in all cases.
In particular, for all block sizes $b$ considered, the effective rank increases less than linearly in the number $d$ of automata and is bounded by $\reff < 20$.
As expected, the maximum rank reacts much more sensitively to changes in $d$ and $b$.

For constant block sizes $b = 2$ and $4$ no smooth increase in the number $d$ of automata can be detected.
Instead, maximum and effective rank oscillate in the number $d$ of automata.
Particularly low ranks typically occur when the strongly correlated automata are not separated by the tree structure.
For example, in case of $d = 32$ automata with a block size of $b = 4$ the strongly correlated automata are split only in the lower levels of the tree, and we observe a very low rank.
For $d = 28$ automata with a block size of $b = 4$, however, the strongly correlated automata $\{\A_{13}, \A_{14}, \A_{15}, \A_{16}\}$, for example, are separated already in the first level of the tree, and we observe a much higher maximum rank.
A closer look at the rank components reveals an increase especially at the points where the blocks are split up.

In the case of $b=\frac{d}{2}$ the block size increases linearly in the number $d$ of automata.
In~\cref{figure:max_rank_b} we see a strong increase of the maximum rank for $d=24, 28, 32$ automata.
However, the effective rank grows smoothly and logarithmically in the number $d$ of automata, see~\cref{figure:eff_rank_b}.
This indicates that the maximum rank occurs only sporadically.
A closer look shows that the maximum ranks again occur only in isolated cases at the vertices where highly correlated automata are separated.
This observation and the fact that the effective rank is small implies that most rank components are small as well.

These results suggest that neither the number $d$ of automata nor the size of the blocks $b$ alone are responsible for an increase in rank, but that especially the distribution of the automata in the tree has a large effect.
How to construct an appropriate tree in order to keep ranks low using a-priori information on the parameters is a topic of ongoing research.
\\

We also  analyze how the approximation rank is influenced by the maximum relative truncation error $\eps$ allowed within the algorithm.
To this end, we compute a low-rank approximation of the marginal distribution $\pPM$ for $d = 32$ automata and blocks of size $b = 4$ using~\cref{alg:uniformtization}. The iteration stops when the norm of the relative residual for~\cref{eq:linear_system} is smaller than $\tol \in \{10^{-1}, 10^{-2}, 10^{-3}, 10^{-4}\}$.
In order to achieve this for all $100$ sample parameters, $\eps$ has to be chosen smaller than $10^{-6}$. If we allow for a higher truncation error, the iteration stagnates in some cases.
We again restrict ourselves to plotting mean values since the variances are small.

\begin{figure}[tbhp]
	\begin{minipage}[tbhp]{0.45\textwidth}
		\vspace{0pt}
		\centering
		
		\includegraphics[width =\textwidth]{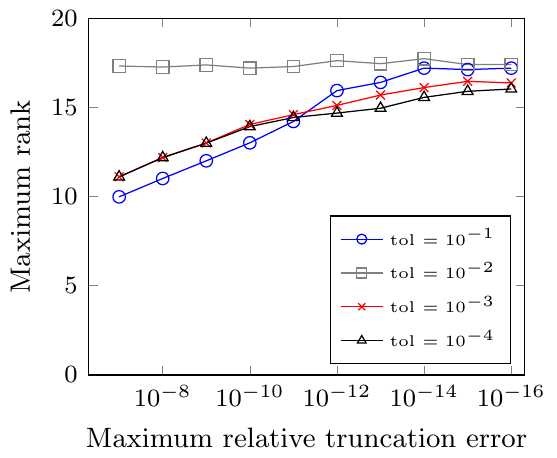}
		
		\subcaption{\footnotesize{Mean of the maximum rank.}}
		\label{figure:max_rank_truncation}
	\end{minipage}
	\hfill
	\begin{minipage}[tbhp]{0.45\textwidth}
		\vspace{0pt}
		\centering

		\includegraphics[width =\textwidth]{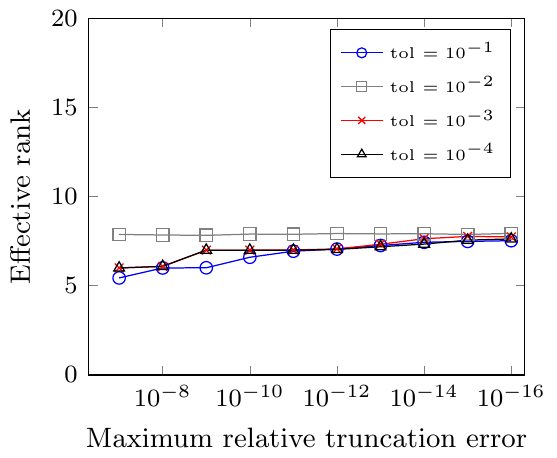}
		
		\subcaption{\footnotesize{Mean of the effective rank.}}
		\label{figure:eff_rank_truncation}
	\end{minipage}
		\caption{\footnotesize{Mean of maximum and effective approximation rank of $\pPM$ as a function of the maximum relative truncation error $\eps$ using~\cref{alg:uniformtization} for $100$ sample parameters with $d=32$ automata and block size $b = 4$. The number of iterations necessary to achieve a relative residual smaller than $\tol$ increases logarithmically in $\tol$ independent of $\eps$.}}
\end{figure}

~\Cref{figure:max_rank_truncation,figure:eff_rank_truncation} show semi-logarithmic plots of the mean of the maximum and effective rank of the approximation of $\pPM$ with $d = 32$ automata and block size $b = 4$ as a function of the maximum relative truncation error $\eps$ for several values of the tolerance $\tol$.
For $\tol = 10^{-2}$ maximum and effective rank are nearly constant with $\rmax \approx 17$ and $\reff \approx 7$.
For $\tol \neq 10^{-2}$, we observe that there are only small differences in maximum and effective rank, respectively.
In these cases, the maximum rank increases logarithmically in the maximum relative truncation error $\eps$, see~\cref{figure:max_rank_truncation}, and the effective rank is nearly constant for all truncation errors, see~\cref{figure:eff_rank_truncation}.
For very small $\eps$ all curves converge.
This observation suggests, on the one hand, that $\pPM$ can be accurately approximated with a tensor of maximum rank $\rmax \approx 17$ and effective rank $\reff \approx 7$, since a more accurate truncation, i.e., smaller $\eps$, has only very small impact on the approximation ranks. 
On the other hand, the conclusion that the ranks are low and nearly independent of the tolerance value indicates that the ranks during the iteration are also low.
This allows not only for efficient storage of the resulting approximation but also for efficient computation using~\cref{alg:uniformtization}.

\subsection{Study of the method}
\label{subsection:study_of_method}

We now study the convergence of~\cref{alg:uniformtization}.

\begin{figure}[tbhp]
	\begin{minipage}[tbhp]{0.45\textwidth}
		\vspace{0pt}
		\centering

		\includegraphics[width = \textwidth]{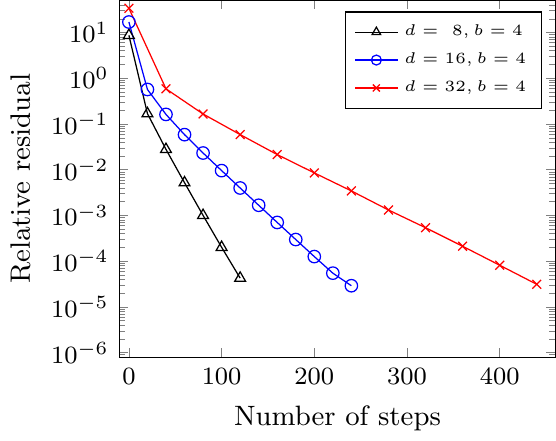}
		
		\subcaption{\footnotesize{Mean for $d \in \{8, 16, 32\}$}}
		\label{figure:convergence_parameter}
	\end{minipage}
	\hfill
	\begin{minipage}[tbhp]{0.45\textwidth}
		\vspace{0pt}
		\centering

		\includegraphics[width = \textwidth]{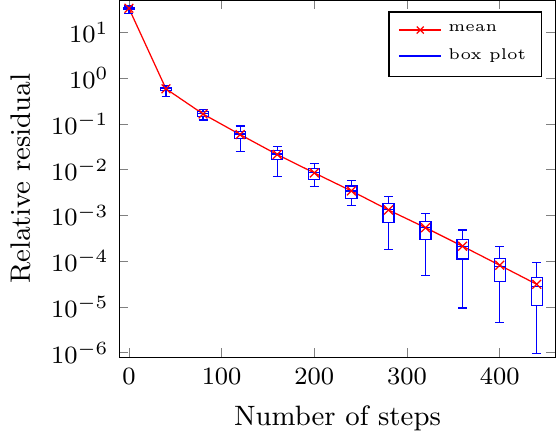}
		
		\subcaption{\footnotesize{Mean and box plot for $d=32$}}
		\label{figure:convergence_rate}
	\end{minipage}
	\caption{\footnotesize{Norm of the relative residual as a function of the iteration step using~\cref{alg:uniformtization} with $\eps = 10^{-8}$ for $100$ sample parameters with $d$ automata and block size $b=4$.}}
\end{figure}

~\Cref{figure:convergence_parameter} and~\cref{figure:convergence_rate} show semi-logarithmic plots of the norm of the relative residual as a function of the iteration step for block size $b = 4$ and $d$ automata. 
~\Cref{figure:convergence_parameter}  displays the mean value of the relative residual and~\cref{figure:convergence_rate} additionally the corresponding box plot illustrating the variances for $d = 32$ automata.
We observe a linear convergence of the method for all values of $d$.
For larger number $d$ of automata the convergence slows down.
This can be explained by the fact that the upper bound $\gamma$ increases in $d$, i.e., the convergence rates $\frac{\gamma}{1 + \gamma}$ are closer to one.
Especially for large $d$ and strongly varying parameters $\PM$ a tighter estimation for the maximum diagonal entries using low-rank formats might be advantageous to speed up the convergence.

In~\cref{figure:convergence_rate} the ranges for $d=32$ given by the boxes are small, which indicates that there are only a few outliers given by the whiskers.
We have confirmed our results using different numbers $d$ of automata and block sizes $b$ (not shown).
The size of the blocks as well as the ranks during the iteration have smaller impact on the convergence.
In our tests we observe that the number of iteration steps needed to achieve a certain tolerance grows linearly in the number $d$ of automata but is independent of the block size $b$.
This indicates that the number of iteration steps is independent of changes in the ordering of automata and consequently of changes in the rank.

\section{Conclusion and future work}
\label{section:conclusion}

Inspired by current research in tumor progression models, we considered a class of continuous-time Markov chains that describe interacting processes and suffer from the problem of state-space explosion.
The corresponding time-marginal distribution is uniquely defined as the solution of a certain linear system.
By representing the Markov chain via a Stochastic Automata Network with separable functional transitions, we obtained a low-rank tensor representation of the operator and the right-hand side of this linear system.
This enabled us to derive an iterative method to compute a low-rank tensor approximation of the time-marginal distribution and hence to overcome the state-space explosion.
 The method guarantees that the entries of the approximation sum up to one as required for a probability distribution.
We proved the convergence of the method.
In numerical experiments focused on the concept of Mutual Hazard Networks we illustrated that the time-marginal distribution is well approximated with low rank.
The method allows for consistently low ranks during the iteration, and linear convergence was observed independently of the number of processes/automata.

A probability distribution, in addition to being normalized to one, must be non-negative. An approximation of a probability distribution should also satisfy these conditions.
How to guarantee non-negativity and at the same time convergence will be part of our future research.
Moreover, we observed that the approximation rank for the time-marginal distribution depends strongly on the structure of the tree tensor network and on the effects between automata.
To minimize the approximation rank we plan to develop a strategy to determine an optimal tree tensor network structure a-priori.

\section*{Acknowledgment}
We thank Tim A. Werthmann for his critical reading of and suggestions for this article.

\bibliographystyle{siamplain}
\bibliography{./tex/references}

\begin{thebibliography}{10}

\bibitem{Anderson2010}
{\sc D.~F. Anderson, G.~Craciun, and T.~G. Kurtz}, {\em {P}roduct-form
  stationary distributions for deficiency zero chemical reaction networks},
  Bulletin of mathematical biology, 72 (2010), pp.~1947--1970,
  \url{https://doi.org/10.1007/s11538-010-9517-4}.

\bibitem{Bachmayr2017}
{\sc M.~Bachmayr and R.~Schneider}, {\em {I}terative {M}ethods {B}ased on
  {S}oft {T}hresholding of {H}ierarchical {T}ensors}, Foundations of
  Computational Mathematics, 17 (2017), pp.~1037--1083,
  \url{https://doi.org/10.1007/s10208-016-9314-z}.

\bibitem{Bailey2018}
{\sc M.~H. e.~a. Bailey}, {\em {C}omprehensive {C}haracterization of {C}ancer
  {D}river {G}enes and {M}utations}, Cell, 173 (2018), pp.~371--385.e18,
  \url{https://doi.org/10.1016/j.cell.2018.02.060}.

\bibitem{Beerenwinkel2014}
{\sc N.~Beerenwinkel, R.~F. Schwarz, M.~Gerstung, and F.~Markowetz}, {\em
  {C}ancer {E}volution: {M}athematical {M}odels and {C}omputational
  {I}nference}, Systematic Biology, 64 (2014), pp.~e1--e25,
  \url{https://doi.org/10.1093/sysbio/syu081}.

\bibitem{Beerenwinkel2009}
{\sc N.~Beerenwinkel and S.~Sullivant}, {\em {{M}arkov models for accumulating
  mutations}}, Biometrika, 96 (2009), pp.~645--661,
  \url{https://doi.org/10.1093/biomet/asp023}.

\bibitem{Kressner2016}
{\sc M.~Bolten, K.~Kahl, D.~Kressner, F.~Macedo, and S.~Sokolović}, {\em
  Multigrid methods combined with low-rank approximation for tensor-structured
  {M}arkov chains}, Electronic Transactions on Numerical Analysis, 48 (2018),
  pp.~348--361, \url{https://doi.org/10.1553/etna_vol48s348}.

\bibitem{Buchholz2007}
{\sc P.~Buchholz and T.~Dayar}, {\em {O}n the {C}onvergence of a {C}lass of
  {M}ultilevel {M}ethods for {L}arge {S}parse {M}arkov {C}hains}, SIAM Journal
  on Matrix Analysis and Applications, 29 (2007), pp.~1025--1049,
  \url{https://doi.org/10.1137/060651161}.

\bibitem{Buchholz2016}
{\sc P.~Buchholz, T.~Dayar, J.~Kriege, and M.~C. Orhan}, {\em {C}ompact
  {R}epresentation of {S}olution {V}ectors in {K}ronecker-{B}ased {M}arkovian
  {A}nalysis}, in Quantitative Evaluation of Systems, G.~Agha and B.~Van~Houdt,
  eds., Cham, 2016, Springer International Publishing, pp.~260--276,
  \url{https://doi.org/10.1007/978-3-319-43425-4_18}.

\bibitem{Buchholz2017}
{\sc P.~Buchholz, T.~Dayar, J.~Kriege, and M.~C. Orhan}, {\em {O}n compact
  solution vectors in {K}ronecker-based {M}arkovian analysis}, Performance
  Evaluation, 115 (2017), pp.~132 -- 149,
  \url{https://doi.org/10.1016/j.peva.2017.08.002}.

\bibitem{Carroll1970}
{\sc J.~D. Carroll and J.-J. Chang}, {\em {A}nalysis of individual differences
  in multidimensional scaling via an n-way generalization of
  ``{E}ckart-{Y}oung'' decomposition}, Psychometrika, 35 (1970), pp.~283--319,
  \url{https://doi.org/10.1007/BF02310791}.

\bibitem{Chan1987}
{\sc R.~H. Chan}, {\em {I}terative methods for overflow queueing models i},
  Numerische Mathematik, 51 (1987), pp.~143--180,
  \url{https://doi.org/10.1007/BF01396747}.

\bibitem{Silva2008}
{\sc V.~de~Silva and L.-H. Lim}, {\em {T}ensor {R}ank and the {I}ll-{P}osedness
  of the {B}est {L}ow-{R}ank {A}pproximation {P}roblem}, SIAM Journal on Matrix
  Analysis and Applications, 30 (2008), pp.~1084--1127,
  \url{https://doi.org/10.1137/06066518X}.

\bibitem{Dubois1979}
{\sc P.~F. Dubois, A.~Greenbaum, and G.~H. Rodrigue}, {\em {A}pproximating the
  inverse of a matrix for use in iterative algorithms on vector processors.},
  Computing, 22 (1979), pp.~257--268, \url{https://doi.org/10.1007/BF02243566}.

\bibitem{Fourneau2008}
{\sc J.-M. Fourneau}, {\em {P}roduct {F}orm {S}teady-{S}tate {D}istribution for
  {S}tochastic {A}utomata {N}etworks with {D}omino {S}ynchronizations}, in
  Proceedings of the 5th European Performance Engineering Workshop on Computer
  Performance Engineering, Berlin, Heidelberg, 2008, Springer-Verlag,
  p.~110–124, \url{https://doi.org/10.1007/978-3-540-87412-6_9}.

\bibitem{Plateau2007}
{\sc J.~M. Fourneau, B.~Plateau, and W.~J. Stewart}, {\em {A}n {A}lgebraic
  {C}ondition for {P}roduct {F}orm in {S}tochastic {A}utomata {N}etworks
  without {S}ynchronizations}, Perform. Eval., 65 (2008), p.~854–868,
  \url{https://doi.org/10.1016/j.peva.2008.04.007}.

\bibitem{Grasedyck2010}
{\sc L.~Grasedyck}, {\em {H}ierarchical {S}ingular {V}alue {D}ecomposition of
  {T}ensors}, SIAM Journal on Matrix Analysis and Applications, 31 (2010),
  pp.~2029--2054, \url{https://doi.org/10.1137/090764189}.

\bibitem{Grasedyck2011}
{\sc L.~Grasedyck and W.~Hackbusch}, {\em {A}n introduction to hierarchical
  ({H}-) rank and {T}{T}-rank of tensors with examples}, Computational Methods
  in Applied Mathematics Comput. Methods Appl. Math., 11 (2011), pp.~291--304,
  \url{https://doi.org/10.2478/cmam-2011-0016}.

\bibitem{Grasedyck2013}
{\sc L.~Grasedyck, D.~Kressner, and C.~Tobler}, {\em {A} literature survey of
  low-rank tensor approximation techniques}, GAMM-Mitteilungen, 36 (2013),
  pp.~53--78, \url{https://doi.org/10.1002/gamm.201310004}.

\bibitem{Grasedyck2017}
{\sc L.~Grasedyck and C.~L{\"o}bbert}, {\em {D}istributed hierarchical {SVD} in
  the {H}ierarchical {T}ucker format}, Numerical Linear Algebra with
  Applications, 25 (2018), \url{https://doi.org/10.1002/nla.2174}.

\bibitem{Grassmann1977}
{\sc W.~Grassmann}, {\em Transient solutions in markovian queueing systems},
  Computers and Operations Research, 4 (1977), pp.~47 -- 53,
  \url{https://doi.org/10.1016/0305-0548(77)90007-7}.

\bibitem{Hackbusch2012}
{\sc W.~Hackbusch}, {\em {T}ensor {S}paces and {N}umerical {T}ensor
  {C}alculus}, vol.~42 of Springer series in computational mathematics,
  Springer, Heidelberg, 2012, \url{https://doi.org/10.1007/978-3-642-28027-6}.

\bibitem{Hackbusch2008}
{\sc W.~Hackbusch, B.~N. Khoromskij, and E.~E. Tyrtyshnikov}, {\em Approximate
  iterations for structured matrices}, Numerische Mathematik, 109 (2008),
  pp.~365--383, \url{https://doi.org/10.1007/s00211-008-0143-0}.

\bibitem{Hackbusch2009}
{\sc W.~Hackbusch and S.~K{\"u}hn}, {\em A {N}ew {S}cheme for the {T}ensor
  {R}epresentation}, Journal of Fourier Analysis and Applications, 15 (2009),
  pp.~706--722, \url{https://doi.org/10.1007/s00041-009-9094-9}.

\bibitem{Lubich2016}
{\sc J.~Haegeman, C.~Lubich, I.~Oseledets, B.~Vandereycken, and F.~Verstraete},
  {\em Unifying time evolution and optimization with matrix product states},
  Phys. Rev. B, 94 (2016), p.~165116,
  \url{https://doi.org/10.1103/PhysRevB.94.165116}.

\bibitem{Harshman1970}
{\sc R.~Harshman}, {\em Foundations of the parafac procedure: models and
  conditions for an 'exploratory' multimodal factor analysis}, in UCLA Working
  Papers in Phonetics, 1970, pp.~1--84.

\bibitem{Hjelm2006}
{\sc M.~Hjelm, M.~H\"{o}glund, and J.~Lagergren}, {\em {N}ew {P}robabilistic
  {N}etwork {M}odels and {A}lgorithms for {O}ncogenesis}, Journal of
  Computational Biology, 13 (2006), pp.~853--865,
  \url{https://doi.org/10.1089/cmb.2006.13.853}.

\bibitem{Hastad1990}
{\sc J.~Håstad}, {\em Tensor rank is {NP}-complete}, Journal of Algorithms, 11
  (1990), pp.~644 -- 654, \url{https://doi.org/10.1016/0196-6774(90)90014-6}.

\bibitem{Johnson2010}
{\sc T.~H. Johnson, S.~R. Clark, and D.~Jaksch}, {\em Dynamical simulations of
  classical stochastic systems using matrix product states}, Physical Review E,
  82 (2010), \url{https://doi.org/10.1103/physreve.82.036702}.

\bibitem{Kazeev2014}
{\sc V.~Kazeev, M.~Khammash, M.~Nip, and C.~Schwab}, {\em {D}irect {S}olution
  of the {C}hemical {M}aster {E}quation {U}sing {Q}uantized {T}ensor {T}rains},
  PLOS Computational Biology, 10 (2014), pp.~1--19,
  \url{https://doi.org/10.1371/journal.pcbi.1003359}.

\bibitem{Kim2014}
{\sc J.~Kim, Y.~He, and H.~Park}, {\em Algorithms for nonnegative matrix and
  tensor factorizations: a unified view based on block coordinate descent
  framework}, Journal of Global Optimization, 58 (2014), pp.~285--319,
  \url{https://doi.org/10.1007/s10898-013-0035-4}.

\bibitem{Kressner2014}
{\sc D.~Kressner and F.~Macedo}, {\em {L}ow-{R}ank {T}ensor {M}ethods for
  {C}ommunicating {M}arkov {P}rocesses}, in Quantitative Evaluation of Systems,
  G.~Norman and W.~Sanders, eds., Cham, 2014, Springer International
  Publishing, pp.~25--40, \url{https://doi.org/10.1007/978-3-319-10696}.

\bibitem{Kressner2012}
{\sc D.~Kressner and C.~Tobler}, {\em Algorithm 941: {H}tucker---{A} {M}atlab
  {T}oolbox for {T}ensors in {H}ierarchical {T}ucker {F}ormat}, ACM Trans.
  Math. Softw., 40 (2014), \url{https://doi.org/10.1145/2538688}.

\bibitem{Kulkarni2011}
{\sc V.~G. Kulkarni}, {\em {I}ntroduction to {M}odeling and {A}nalysis of
  {S}tochastic {S}ystems}, Springer-Verlag New York, 2~ed., 2011,
  \url{https://doi.org/10.1007/978-1-4419-1772-0}.

\bibitem{Masetti2019}
{\sc G.~Masetti and L.~Robol}, {\em {T}ensor methods for the computation of
  {MTTF} in large systems of loosely interconnected components}, tech. report,
  ISTI-CNR Open Portal, 2019, \url{https://arxiv.org/abs/1907.02449}.

\bibitem{Matlab2019}
{\sc {Natick, Massachusetts}}, {\em {MATLAB}}, 2019,
  \url{https://de.mathworks.com/products/matlab.html} (accessed 2020-05-05).
\newblock Version 9.6 (R2019a).

\bibitem{Oseledets2009}
{\sc I.~V. Oseledets and E.~Tyrtyshnikov}, {\em {B}reaking the {C}urse of
  {D}imensionality, {O}r {H}ow to {U}se {SVD} in {M}any {D}imensions}, SIAM
  Journal on Scientific Computing, 31 (2009), pp.~3744--3759,
  \url{https://doi.org/10.1137/090748330}.

\bibitem{Ostlund1995}
{\sc S.~\"Ostlund and S.~Rommer}, {\em Thermodynamic limit of density matrix
  renormalization}, Phys. Rev. Lett., 75 (1995), pp.~3537--3540,
  \url{https://doi.org/10.1103/PhysRevLett.75.3537}.

\bibitem{Plateau2000}
{\sc B.~Plateau and W.~J. Stewart}, {\em {S}tochastic {A}utomata {N}etworks},
  Springer US, Boston, MA, 2000, pp.~113--151,
  \url{https://doi.org/10.1007/978-1-4757-4828-4_5}.

\bibitem{Schill2019}
{\sc R.~Schill, S.~Solbrig, T.~Wettig, and R.~Spang}, {\em {{M}odelling cancer
  progression using {M}utual {H}azard {N}tworks}}, Bioinformatics, 36 (2019),
  pp.~241--249, \url{https://doi.org/10.1093/bioinformatics/btz513}.

\bibitem{Loan2008}
{\sc C.~F. Van~Loan}, {\em {T}ensor {N}etwork {C}omputations in {Q}uantum
  {C}hemistry}, 2008, \url{www.cs.cornell.edu/cv/OtherPdf/ZeuthenCVL.pdf}
  (accessed 2020-05-05).

\bibitem{White1992}
{\sc S.~R. White}, {\em Density matrix formulation for quantum renormalization
  groups}, Phys. Rev. Lett., 69 (1992), pp.~2863--2866,
  \url{https://doi.org/10.1103/PhysRevLett.69.2863}.

\bibitem{Yamanaka1997}
{\sc K.~Yamanaka, M.~Agu, and T.~Miyajima}, {\em {A} {C}ontinuous-{T}ime
  {A}synchronous {B}oltzmann {M}achine}, Neural Networks, 10 (1997),
  pp.~1103--1107, \url{https://doi.org/10.1016/s0893-6080(97)00006-3}.

\end{thebibliography}

\end{document}